\documentclass[11pt]{article}
  \usepackage{array}                        
  \usepackage{theorem}                      
  \usepackage{amsmath,amscd,amssymb}        
  \usepackage{latexsym}                     
  \usepackage[german,english]{babel}        
  \usepackage[applemac]{inputenc}           
\usepackage{xypic}

\newtheorem{theorem}{Theorem}[section]
\newtheorem{corollary}[theorem]{Corollary}
\newtheorem{proposition}[theorem]{Proposition}
\newtheorem{lemma}[theorem]{Lemma}
\theorembodyfont{\rmfamily}

\newtheorem{rem}[theorem]{Remark}
\newtheorem{urems}{Remarks}  
\newtheorem{conjecture}{Conjecture}  
\newtheorem{acknowledgement}{Acknowledgement}

\newcommand{\CC}{{\mathbb{C}}}
\newcommand{\PP}{{\mathbb{P}}}
\newcommand{\QQ}{{\mathbb{Q}}}

\newcommand{\HH}{{\mathbb{H}}}
\newcommand{\RR}{{\mathbb{R}}}
\newcommand{\ZZ}{{\mathbb{Z}}}

\newcommand{\calA}{{\cal A}}

\newcommand{\calM}{{\cal M}}

\newcommand{\calD}{{\cal D}}

\newcommand{\eps}{\varepsilon}

\newcommand{\on}[1]{\operatorname{#1}}
\newenvironment{Proof}{\begin{ProofwCaption}{Proof}}{\end{ProofwCaption}}
\newenvironment{Proof*}[1]{\begin{ProofwCaption}{{#1}}}{\end{ProofwCaption}}
\newenvironment{ProofwCaption}[1]%
  {\addvspace\theorempreskipamount \noindent{\it #1.}\rm}%
  {\qed \par \addvspace\theorempostskipamount}
\newcommand{\qedsymbol}{\mbox{$\Box$}}
\newcommand{\qed}{\hfill\qedsymbol}

\begin{document}
\title{Nef Divisors on Moduli Spaces of Abelian Varieties}
\author{K.~Hulek}
\date{\em{Dedicated to the memory of Michael Schneider}}
\maketitle

%
%
\setcounter{section}{-1}
\section{Introduction}

Let $\calA_g$ be the moduli space of principally polarized abelian varieties
of dimension $g$. Over the complex numbers $\calA_g=\HH_g/\Gamma_g$ where
$\HH_g$ is the Siegel space of genus $g$ and $\Gamma_g=\on{Sp}(2g,\ZZ)$. We
denote the torodial compactification given by the second Voronoi decomposition
by $\calA_g^*$ and call it the \emph{Voronoi compactification}.
It was shown by
Alexeev and Nakamura \cite{A} that $\calA_g^*$ coarsely represents the stack
of principally polarized stable quasiabelian varieties. The variety $\calA_g^*$
is projective \cite{A} and it is known that
the Picard group of ${\cal A}_g^*, g\ge 2$ is generated
(modulo
torsion) by two elements
$L$ and
$D$, where $L$ denotes the ($\QQ$-)line bundle given by modular forms of
weight $1$ and
$D$ is the boundary (see \cite{Mu2}, \cite{Fa} and \cite{Mu1} for $g=2, 3$ and
$\geq4$).
In this paper we want to discuss the
following

\begin{theorem}\label{theo0.1}
Let $g=2$ or $3$. A divisor $aL-bD$ on $\calA_g^*$ is nef if and only if
$b\geq0$ and $a-12b\geq0$.
\end{theorem}

The varieties $\calA_g$ have finite quotient singularities. Adding a level-$n$
structure one obtains spaces $\calA_g(n)=\HH_g/\Gamma_g(n)$ where
$\Gamma_g(n)$ is the principal congruence subgroup of level $n$. For $n\geq3$
these spaces are smooth. However, the Voronoi compactification $\calA_g^*(n)$
acquires singularities on the boundary for $g\geq5$ due to bad behaviour of the
second Voronoi decomposition. There is a natural quotient map
$\calA_g^*(n)\to\calA_g^*$. Note that this map is branched of order $n$ along
the boundary. Hence Theorem~(\ref{theo0.1}) is equivalent to

\begin{theorem}\label{theo0.2}
Let $g=2$ or $3$. A divisor $aL-bD$ on $\calA_g^*(n)$ is nef if and only if
$b\geq0$ and $a-12\frac{b}{n}\geq0$.
\end{theorem}

This theorem easily gives the following two corollaries.

\begin{corollary}\label{cor0.3}
If $g=2$ then $K$ is nef but not ample for $\calA_2^*(4)$ and $K$ is ample for
$\calA_2^*(n)$, $n\geq5$; in particular $\calA_2^*(n)$ is a minimal model for
$n\geq4$ and a canonical model for $n\geq5$.
\end{corollary}

This was first proved by Borisov \cite{Bo}.

\begin{corollary}\label{cor0.4}
If $g=3$ then $K$ is nef but not ample for $\calA_3^*(3)$ and $K$ is ample for
$\calA_3^*(n)$, $n\geq4$; in particular $\calA_3^*(n)$ is a minimal model for
$n\geq3$ and a canonical model for $n\geq4$.
\end{corollary}

In this paper we shall give two proofs of Theorem~(\ref{theo0.1}). The first
and quick one reduces the problem via the Torelli map to the analogous
question for ${\overline{M}}_2$, resp.~${\overline{M}}_3$. Since the
Torelli map is not surjective for $g\ge 4$ this proof cannot possibly be
generalized to higher genus. This is the main reason why we want to give a
second
proof which uses theta functions. This proof makes essential use of a result
of Weissauer
\cite{We}. The method has the advantage that it extends in principle to other
polarizations as well as to higher
$g$.
We will also give some partial results supporting the

\begin{conjecture}
For any $g\geq2$ the nef cone on $\calA_g^*$ is given by the divisors $aL-bD$
where $b\geq0$ and $a-12b\geq0$.
\end{conjecture}

\begin{acknowledgement}
It is a pleasure for me to thank RIMS and Kyoto University for their
hospitality during the autumn of 1996. I am grateful to V.~Alexeev and
R.~Salvati Manni for useful discussions. It was Salvati Manni who drew my
attention to Weissauer's paper.
I would also like to thank
R.~Weissauer for additional information on
\cite{We}. The author is partially supported by TMR
grant ERBCHRXCT 940557.
\end{acknowledgement}

%
%
\section{Curves meeting the interior}

We start by recalling some results about the Kodaira dimension of
$\calA_g^*(n)$. It was proved by Freitag, Tai and Mumford that $\calA_g^*$ is
of general type for $g\geq7$. The following more general result is probably
well known to some specialists.

\begin{theorem}\label{theo1.1}
$\calA_g^*(n)$ is of general type for the following values of $g$ and $n\geq
n_0$:
\begin{center}
\begin{tabular}{c|cccccc}
$g$   & $2$ & $3$ & $4$ & $5$ & $6$ & $\geq7$ \\\hline
$n_0$ & $4$ & $3$ & $2$ & $2$ & $2$ & $1$
\end{tabular}.
\end{center}
\end{theorem}

\begin{Proof}
One can use Mumford's method from \cite{Mu1}. First recall that away from the
singularities and the closure of the branch locus of the
map $\HH_g\to\calA_g(n)$ the canonical
bundle equals
\begin{equation}\label{formula1.1}
K\equiv (g+1)L-D.
\end{equation}
This equality holds in particular also on an open part of the boundary.
If $g\leq4$ and
$n\geq3$ the spaces $\calA_g^*(n)$ are smooth and hence (\ref{formula1.1})
holds
everywhere. If $g\geq5$ then Tai \cite{T} showed that there is a suitable
toroidal compactification ${\tilde{\calA}}_g(n)$ such that all singularities
are canonical quotient singularities. By Mumford's results from \cite{Mu1} one
can use the theta-null locus to eliminate $D$ from formula (\ref{formula1.1})
and obtains
\begin{equation}\label{formula1.2}
K\equiv \left((g+1)-\frac{2^{g-2}(2^g+1)}{n2^{2g-5}}\right)L+
\frac{1}{n2^{2g-5}}[\Theta_{\on{null}}].
\end{equation}
We then have general type if all singularities are canonical and if the factor
in front of $L$ is positive. This gives immediately all values in the above
table with the exception of $(g,n)=(4,2)$ and $(7,1)$. In the latter case
the factor in front of $L$ is negative. The proof that ${\cal A}_7$
is nevertheless of general type is the
main result of \cite{Mu1}. The difficulty in the first case is
that one can possibly have non-canonical singularities.
One can, however, use the following argument which I
have learnt from Salvati Manni: An immediate calculation shows that for every
element $\sigma\in\Gamma_g(2)$ the square $\sigma^2\in\Gamma_g(4)$. Hence if
$\sigma$ has a fixed point then $\sigma^2=1$ since $\Gamma_g(4)$ acts freely.
But for elements of order $2$ one can again use Tai's extension theorem (see
\cite[Remark after Lemma~4.5]{T} and \cite[Remark after Lemma~5.2]{T}).
\end{Proof}

\begin{rem}\label{rem1.2}
The Kodaira dimension of $\calA_6$ is still unknown. All other
varieties ${\cal A}_g(n)$ which do not appear in the above list are either
rational or unirational: Unirationality of
$\calA_g$ for $g=5$ was proved  by Donagi \cite{D} and by Mori and Mukai
\cite{MM}. For $g=4$ the same result was shown by Clemens \cite{C}.
Unirationality is
easy for
$g\leq3$. Igusa
\cite{I2} showed that $\calA_2$ is rational. Recently Katsylo \cite{Ka} proved
rationality of ${\cal M}_3$ and hence also of ${\cal A}_3$.
The space
$\calA_3(2)$ is rational by work of van Geemen \cite{vG} and
Dolgachev and Ortlang \cite {DO}. $\calA_2(3)$ is the Burkhardt quartic
and hence
rational. This was first
proved by Todd
(1936)  and Baker (1942). See also the thesis of Finkelnberg
\cite{Fi}.
The variety ${\cal A}_2(2)$ has the Segre cubic as a projective
model \cite{vdG1} and is hence also rational.
Yamazaki
\cite{Ya} first showed general type for $\calA_2(n)$, $n\geq4$.
\end{rem}

We denote the Satake compactification of $\calA_g$ by $\overline{\calA}_g$.
There is a natural map $\pi:\calA_g^*\to\overline{\calA}_g$ which is an
isomorphism on $\calA_g$. The line bundle $L$ is the pullback of an ample line
bundle on $\overline{\calA}_g$ which, by abuse of notation, we again denote by
$L$. In fact the Satake compactification is defined as the closure of
the image
of $\calA_g$ under the embedding given by a suitable power of $L$ on $\calA_g$.
In particular we notice that $L.C\geq0$ for every curve $C$ on $\calA_g^*$ and
that $L.C>0$ if $C$ is not contracted to a point under the map $\pi$.

Let $F$ be a modular form with respect to the full modular group
$\on{Sp}(2g,\ZZ)$. Then the \emph{order} $o(F)$ of $F$ is defined as the
quotient of the vanishing order of $F$ divided by the weight of $F$.

\begin{theorem}[Weissauer]\label{theo1.3}
For every point $\tau\in\HH_g$ and every $\eps>0$ there exists a modular form
$F$ of order $o(F)\geq\frac{1}{12+\eps}$ which does not vanish at $\tau$.
\end{theorem}

\begin{Proof}
See \cite{We}.
\end{Proof}

\begin{proposition}\label{prop1.4}
Let $C\subset\calA_g^*$ be a curve which is not contained in the boundary.
Then $(aL-bD).C\geq0$ if $b\geq0$ and $a-12b\geq0$.
\end{proposition}

\begin{Proof}
First note that $L.C>0$ since $\pi(C)$ is a curve in the Satake
compactification. It is enough to prove that $(aL-bD).C>0$ if $a-12b>0$ and
$a,b\geq0$. This is clear for $b=0$ and hence we can assume that $b\neq0$. We
can now choose some $\eps>0$ with $a/b>12+\eps$. By Weissauer's theorem there
exists a modular form $F$ of say weight $k$ and vanishing order $m$ with
$F(\tau)\neq0$ for some point $[\tau]\in C$ and $m/k\geq1/(12+\eps)$. In terms
of divisors this gives us that
$$
kL=mD+D_F, \quad C\not\subset D_F
$$
where $D_F$ is the zero-divisor of $F$. Hence
$$
\left(\frac km L-D\right)=\frac1mD_F.C\geq0.
$$
Since $a/b>12+\eps\geq k/m$ and $L.C>0$ we can now conclude that
$$
\left(\frac abL-D\right).C>\left(\frac kmL-D\right).C\geq0.
$$
\end{Proof}

\begin{rem}\label{rem1.5}
Weissauer's result is optimal, since the modular forms of order $>1/12$
have a common base locus. To see this consider curves $C$ in
$\calA_g^*$ of the form $X(1)\times\{A\}$ where $X(1)$ is the modular curve of
level $1$ parametrizing elliptic curves and $A$ is a fixed abelian variety of
dimension $g-1$. The degree of $L$ on $X(1)$ is $1/12$ (recall that $L$ is a
$\QQ$-bundle) whereas it has one cusp, i.e.~the degree of $D$ on this curve is
$1$. Hence every modular form of order $>1/12$ will vanish on $C$. This also
shows that the condition $a-12b\geq0$ is necessary for a
divisor to be nef.
\end{rem}

%
%
\section{Geometry of the boundary (I)}
\setcounter{equation}{0}

We first have to collect some properties of the structure of the boundary of
$\calA_g^*(n)$. Recall that the Satake compactification is set-theoretically
the union of $\calA_g(n)$ and of moduli spaces $\calA_k(n)$, $k<g$ of lower
dimension, i.e.
$$
\overline{\calA}_g(n)=\calA_g(n)\amalg\left(\underset{i_1}{\amalg}
\calA_{g-1}^{i_1}(n)\right)\amalg \left(\underset{i_2}{\amalg}
\calA_{g-2}^{i_2}(n)\right)\ldots\amalg \left(\underset{i_g}{\amalg}
\calA_0^{i_g}(n)\right).
$$
Via the map $\pi:\calA_g^*(n)\to\overline{\calA}_g(n)$ this also defines a
stratification of $\calA_g^*(n)$:
$$
\calA^*_g(n)=\calA_g(n)\amalg\left(\underset{i_1}{\amalg}
D_{g-1}^{i_1}(n)\right)\amalg \left(\underset{i_2}{\amalg}
D_{g-2}^{i_2}(n)\right)\ldots\amalg \left(\underset{i_g}{\amalg}
D_0^{i_g}(n)\right).
$$
The irreducible components of the boundary $D$ are the closures
$\overline{D}_{g-1}^{i_1}(n)$ of the codimension $1$
strata $D_{g-1}^{i_1}(n)$. Whenever we
talk about a \emph{boundary component} we mean one of
the divisors $\overline{D}_{g-1}^{i_1}(n)$. Then the boundary $D$ is given by
$$
D=\sum_{i_1}\overline{D}_{g-1}^{i_1}(n).
$$
The fibration $\pi:D_{g-1}^{i_1}(n)\to\calA_{g-1}^{i_1}(n)=\calA_{g-1}(n)$ is
the universal family of abelian varieties of dimension $g-1$ with a level-$n$
structure if $n\geq3$ resp.~the universal family of Kummer surfaces for $n=1$
or $2$ (see \cite{Mu1}). We shall also explain this in more detail later on.
To be more precise we associate to a point $\tau\in\HH_g$ the lattice
$L_{\tau,1}=(\tau,\mathbf{1})\ZZ^{2g}$, resp.~the principally polarized
abelian variety $A_{\tau,1}=\CC^g/L_{\tau,1}$. Given an integer $n\geq1$ we
set $L_{n\tau,n}=(n\tau,n\mathbf{1}_g)\ZZ^{2g}$, resp.~$A_{n\tau,n}=
\CC^g/L_{n\tau,n}$. By $K_{n\tau,n}$ we denote the Kummer variety $A_{n,\tau
n}/\{\pm1\}$.

\begin{lemma}\label{lemma2.1}
Let $n\geq3$. Then for any point $[\tau]\in\calA_{g-1}^{i_1}(n)$ the fibre of
$\pi$ equals $\pi^{-1}([\tau])=A_{n,\tau n}$.
\end{lemma}

\begin{Proof}
Compare \cite{Mu1}. We shall also give an independent proof below.
\end{Proof}

This result remains true for $n=1$ or $2$, at least for points $\tau$ whose
stabilizer subgroup in $\Gamma_g(n)$ is $\{\pm1\}$, if we replace $A_{n,\tau
n}$ by its associate Kummer variety $K_{n,\tau n}$.

\begin{lemma}\label{lemma2.2}
Let $n\geq3$. Then for $[\tau]\in\calA_{g-1}^{i_1}(n)$ the restriction of
${D}_{g-1}^{i_1}(n)$
to the fibre $\pi^{-1}([\tau])$ is negative. More precisely
$$
D_{g-1}^{i_1}(n)|_{\pi^{-1}([\tau])}\equiv -\frac2nH
$$
where $H$ is the polarization on $A_{n\tau,n}$ given by the pull-back of the
principal polarization on $A_{\tau,1}$ via the covering $A_{n,\tau n}\to
A_{\tau,1}$.
\end{lemma}

\begin{Proof}
Compare \cite[Proposition~1.8]{Mu1}, resp.~see the discussion below.
\end{Proof}

Again the statement remains true for $n=1$ or $2$ if we replace the abelian
variety by its Kummer variety.

\begin{Proof*}{First proof of Theorem (\ref{theo0.1})}
We have already seen (see Remark~\ref{rem1.5}) that for
every nef divisor $aL-bD$ the inequality
$a-12b\geq0$ holds. If $C$ is a curve in a fibre of the map
$\calA_g^*(n)\to\overline{\calA}_g(n)$, then $L.C=0$. Lemma~(\ref{lemma2.2})
immediately implies that $b\geq0$ for any nef divisor. It remains to show that
the conditions of Theorem~(\ref{theo0.1}) are sufficient to imply nefness.
For any genus the Torelli map $t:\calM_g\to\calA_g$ extends to a morphism
$\overline{t}:\overline{\calM}_g\to \calA_g^*$ (see \cite{Nam}).
Here $\overline{\calM}_g$ denotes the compactification of $\calM_g$ by stable
curves. For $g=2$ and $3$ the map $\overline{t}$ is surjective.
It follows that for every curve $C$ in $\calA_g^*$ there is a curve $C'$
in $\overline{\calM}_g$ which is finite over $C$.
Hence a divisor on $\calA_g^*$, $g=2,3$ is nef if and only if this
holds for its pull-back to $\overline{\calM}_g$.
In the notation of Faber's paper
\cite{Fa} $\overline{t}^*L=\lambda$ where $\lambda$ is the Hodge bundle and
$\overline{t}^*D=\delta_0$ where
$\delta_0$ is the boundary ($g=2$), resp. the closure of the locus of
genus $2$ curves with one node ($g=3$) (cf also \cite{vdG2}). The
result follows since $a\lambda-b\delta_0$ is nef on $\overline{\calM}_g$,
$g=2,3$ for $a-12b\geq0$ and $b\geq0$ (see \cite{Fa}).
\end{Proof*}

As we have already pointed out the Torelli map is not surjective for
$g\geq4$ and hence this proof cannot possibly be generalized to
higher genus.
The main purpose of this paper is, therefore, to give a proof of
Theorem~(\ref{theo0.1}) which does not use the reduction to the curve case.
This
will also allow us to prove some results for general $g$. At the same time we
obtain an independent proof of nefness of $a\lambda-b\delta_0$ for
$a-12b\geq0$ and $b\geq0$ on $\overline{\calM}_g$ for $g=2$ and $3$.

We now want to investigate the open parts $D_{g-1}^{i_1}(n)$ of the boundary
components $\overline{D}_{g-1}^{i_1}(n)$ and their fibration over
$\calA_{g-1}(n)$ more closely. At the same time this gives us another
argument for
Lemmas~(\ref{lemma2.1}) and (\ref{lemma2.2}). At this stage we have to make
first use of the toroidal construction. Recall that the boundary components
$D_{g-1}^{i_1}(n)$  are in $1:1$ correspondence with the maximal dimensional
cusps, and these in turn are in $1:1$ correspondence with the lines
$l\subset\QQ^g$ modulo $\Gamma_g(n)$. Since all cusps are equivalent under the
action of $\Gamma_g/\Gamma_g(n)$ we can restrict our attention to one of these
cusps, namely the one given by $l_0=(0,\ldots,0,1)$. This corresponds to
$\tau_{gg}\to i\infty$. To simplify notation we shall denote the corresponding
boundary stratum simply by ${D}_{g-1}^{1}(n)=D_{g-1}(n)$.
The stabilizer $P(l_0)$ of $l_0$ in
$\Gamma_g$ is generated by elements of the following form
(cf.~\cite[Proposition~I.3.87]{HKW}):
\begin{align*}
g_1&=\begin{pmatrix} A&0&B&0 \\ 0&1&0&0 \\ C&0&D&0 \\ 0&0&0&1 \end{pmatrix},\
\begin{pmatrix} A&B \\ C&D \end{pmatrix} \in \Gamma_{g-1},\\
g_2&=\begin{pmatrix} \mathbf{1}_{g-1}&0&0&0 \\ 0&\pm1&0&0 \\
0&0&\mathbf{1}_{g-1}&0 \\ 0&0&0&\pm1 \end{pmatrix},\\
g_3&=\begin{pmatrix} \mathbf{1}_{g-1}&0&0&\sideset{^t}{}{\on{\mathit{N}}} \\
M&1&N&0 \\ 0&0&\mathbf{1}_{g-1}&-\sideset{^t}{}{\on{\mathit{M}}} \\ 0&0&0&1
\end{pmatrix},\ M,N \in\ZZ^{g-1},\\
g_4&=\begin{pmatrix} \mathbf{1}_{g-1}&0&0&0 \\ 0&1&0&S \\
0&0&\mathbf{1}_{g-1}&0 \\ 0&0&0&1 \end{pmatrix},\ S\in\ZZ.
\end{align*}

\noindent We write $\tau=(\tau_{ij})_{1\leq i,j\leq g}$ in the form
$$
\left(
\begin{array}{ccc|c}
\tau_{11}    & \cdots & \tau_{1,g-1}   & \tau_{1g}\\
\vdots       &        & \vdots         & \vdots\\
\tau_{1,g-1} & \cdots & \tau_{g-1,g-1} & \tau_{g-1,g}\\\hline
\tau_{1,g}   & \cdots & \tau_{g-1,g}   & \tau_{gg}
\end{array}
\right)
=
\left(
\begin{array}{c|c}
\tau_1 & \sideset{^t}{}{\on{\tau_2}}\\\hline
\tau_2 & \tau_3
\end{array}
\right).
$$

\noindent Then the action of $P(l_0)$ on $\HH_g$ is given by
(cf.~\cite[I.3.91]{HKW}):
\begin{align*}
g_1(\tau)&=\begin{pmatrix} (A\tau_1+B)(C\tau_1+D)^{-1} & * \\
\tau_2(C\tau_1+D)^{-1} & \tau_3-\tau_2(C\tau_1+D)^{-1}C
\sideset{^t}{}{\on{\tau_2}} \end{pmatrix},\\
g_2(\tau)&=\begin{pmatrix} \tau_1 & * \\ \pm\tau_2 & \tau_3 \end{pmatrix},\\
g_3(\tau)&=\begin{pmatrix} \tau_1 & * \\ \tau_2+M\tau_1+N & \tau_3'
\end{pmatrix}\\
\intertext{where $\tau_3'=\tau_3+M\tau_1\sideset{^t}{}{\on{\mathit{M}}} +
M\sideset{^t}{}{\on{\tau_2}}+ \sideset{^t}{}{\on{(}} M
\sideset{^t}{}{\on{\tau_2}})+N\sideset{^t}{}{\on{\mathit{M}}},$}
g_4(\tau)&=\begin{pmatrix} \tau_1 & \tau_2
\\ \tau_2 &
\tau_3+S
\end{pmatrix}.
\end{align*}

The parabolic subgroup $P(l_0)$ is an extension
$$
1\longrightarrow P'(l_0)\longrightarrow P(l_0)\longrightarrow P''(l_0)
\longrightarrow 1
$$
where $P'(l_0)$ is the rank $1$ lattice generated by $g_4$. To obtain the same
result for $\Gamma_g(n)$ we just have to intersect $P(l_0)$ with
$\Gamma_g(n)$. Note that $g_2$ is in $\Gamma_g(n)$ only for $n=1$ or $2$. The
first step in the construction of the toroidal compactification
of $\calA_g^*(n)$
is to divide $\HH_g$ by $P'(l_0)\cap\Gamma(n)$ which gives a map
$$
\begin{array}{ccl}
\HH_g &\longrightarrow& \HH_{g-1}\times\CC^{g-1}\times\CC^*\\
\begin{pmatrix} \tau_1 & \sideset{^t}{}{\on{\tau_2}} \\ \tau_2 & \tau_3
\end{pmatrix} &\longmapsto& (\tau_1,\tau_2,e^{2\pi i\tau_3/n}).
\end{array}
$$
\noindent Partial compactification in the direction of $l_0$ then consists of
adding the set $\HH_{g-1}\times\CC^{g-1}\times\{0\}$. It now follows
immediately
from the above formulae for the action of $P(l_0)$ on $\HH_g$ that
the action of
the  quotient group $P''(l_0)$ on  $\HH_{g-1}\times\CC^{g-1}\times\CC^*$
extends
to
$\HH_{g-1}\times\CC^{g-1}\times\{0\}$.
Then
$D_{g-1}(n)=(\HH_{g-1}\times\CC^{g-1})/P''(l_0)$ and the map to
$\calA_{g-1}(n)$ is induced by the projection from $\HH_{g-1}\times\CC^{g-1}$
to $\HH_{g-1}$. This also shows that $D_{g-1}(n)\to\calA_{g-1}(n)$ is the
universal family for $n\geq3$ and that the general fibre is a Kummer variety
for $n=1$ and $2$.

Whenever $n_1|n_2$ we have a Galois covering
$$
\pi(n_1,n_2):\calA_g^*(n_2)\longrightarrow\calA_g^*(n_1)
$$
whose Galois group is $\Gamma_g(n_1)/\Gamma_g(n_2)$. This induces coverings
$\overline{D}_{g-1}(n_2)\to\overline{D}_{g-1}(n_1)$, resp.~$D_{g-1}(n_2)\to
D_{g-1}(n_1)$. In order to avoid technical difficulties we assume for
the moment that $\calA_g^*(n)$ is smooth (this is the case if $g\leq4$
and $n\geq3$). In what follows we will always be able to assume that we are
in this situation.
Then we denote the normal bundle of $\overline{D}_{g-1}(n)$ in
$\calA_g^*(n)$ by $N_{\overline{D}_{g-1}(n)}$, resp.~its restriction to
$D_{g-1}(n)$ by $N_{D_{g-1}(n)}$. Since the covering map $\pi(n_1,n_2)$ is
branched of order $n_2/n_1$ along the boundary, it follows that
$$
\pi^*(n_1,n_2)n_1N_{\overline{D}_{g-1}}(n_1)=n_2N_{\overline{D}_{g-1}}(n_2).
$$
We now define the bundle
$$
\overline{M}(n):=-nN_{\overline{D}_{g-1}(n)}+L.
$$
This is a line bundle on the boundary component $\overline{D}_{g-1}(n)$. We
denote the restriction of $\overline{M}(n)$ to $D_{g-1}(n)$ by $M(n)$. We find
immediately that
$$
\pi^*(n_1,n_2)\overline{M}(n_1)=\overline{M}(n_2).
$$

The advantage of working with the bundle $\overline{M}(n)$ is that we can
explicitly describe sections of this bundle.
For this purpose it is useful to review some basic facts about theta functions.
For every element $m=(m',m'')$ of $\RR^{2g}$ one can define the theta-function
$$
\Theta_{m'm''}(\tau,z)=\sum_{q\in\ZZ^g} e^{2\pi i[{(q+m')\tau}
{^t(q+m')}/2+{(q+m')}{^t(z+m'')}]}.
$$
The transformation behaviour of $\Theta_{m'm''}(\tau,z)$ with respect to
$z\mapsto z+u\tau+u'$ is described by the formulae ($\Theta1$)--($\Theta5$) of
\cite[pp.~49,~50]{I}. The behaviour of $\Theta_{m'm''}(\tau,z)$
with respect to
the action of $\Gamma_g(1)$ on $\HH_g\times\CC^g$ is given by the theta
transformation formula \cite[Theorem~II.5.6]{I} resp.~the corollary following
this theorem \cite[p.~85]{I}.

\begin{proposition}\label{prop2.3}
Let $n\equiv0\on{mod}4p^2$. If $m',m'',\overline{m}',\overline{m}''\in
\frac{1}{2p} \ZZ^{g-1}$, then the functions $\Theta_{m'm''}(\tau,z)
\Theta_{\overline{m}'\overline{m}''}(\tau,z)$ define sections of the line
bundle $M(n)$ on $D_{g-1}(n)$.
\end{proposition}

\begin{Proof}
It follows from ($\Theta3$) and ($\Theta1$) that for $k,k'\in n\ZZ^{g-1}$ the
following holds:
$$
\Theta_{m',m''}(\tau,z+k\tau+k')= e^{2\pi i[-\frac12k{\tau}
{^tk}-{k}{^t(z+k')}]} \Theta_{m',m''}(\tau,z).
$$
Similarly, of course, for $\Theta_{\overline{m}',\overline{m}''}(\tau,z)$.
Moreover the theta transformation formula together with formula ($\Theta2$)
gives
$$
\Theta_{m',m''}(\tau^\#,z^\#) =
e^{2\pi i[\frac12z(C\tau+D)^{-1}{C}{^tz}]}
\det(C\tau+D)^{1/2} u \Theta_{m',m''}(\tau,z)
$$
for every element $\gamma=\begin{pmatrix} A&B\\ C&D \end{pmatrix}\in
\Gamma_{g-1}(n)$ and
$$
\tau^\#=\gamma(\tau),\quad z^\#=z(C\tau+D)^{-1}.
$$
Here $u^2$ is a character of
$\Gamma_{g-1}(1,2)$ with $u^2|_{\Gamma_{g-1}(4)}\equiv1$.

On the other hand the boundary component $D_{g-1}(n)$ is defined by $t_3=0$
with $t_3=e^{2\pi i\tau_3/n}$. We have already described the action of
$P''(l_0)$ on $\HH_{g-1}\times\CC^{g-1}$. The result then follows by comparing
the transformation behaviour of $(t_3/t_3^2)^n$
with respect to $g_1$ and $g_3$
with the above formulae together with the fact that the line bundle $L$ is
defined by the automorphy factor $\det(C\tau+D)$.
\end{Proof}

This also gives an independent proof of Lemma~(\ref{lemma2.2}).

\section{Geometry of the boundary (II)}

So far we have described the stratum $D_{g-1}(n)$ of the boundary component
$\overline{D}_{g-1}(n)$ and we have seen that there is a natural map
$D_{g-1}(n)\to\calA_{g-1}(n)$ which identifies $D_{g-1}(n)$ with the universal
family over $\calA_{g-1}(n)$ if $n\geq3$. We now want to describe the closure
$\overline{D}_{g-1}(n)$ in some detail. In order to do this we have to
restrict ourselves to $g=2$ and $3$. First assume $g=2$. Then the projection
$D_1(n)\to\calA_1(n)=X^0(n)$ extends to a projection
$\overline{D}_1(n)\to X(n)$
onto the modular curve of level $n$ and in this way $\overline{D}_1(n)$ is
identified with Shioda's modular surface $S(n)\to X(n)$. The fibres are either
elliptic curves or $n$-gons of rational curves (if $n\ge 3$).
Similarly the fibration
$D_2(n)\to\calA_2(n)$ extends to a fibration $\overline{D}_2(n)\to\calA_2^*$
whose fibres over the boundary of $\calA_2^*(n)$ are degenerate abelian
surfaces. This was first observed by Nakamura \cite{Na} and was described in
detail by Tsushima \cite{Ts} whose paper is essential for what follows.

We shall now explain the toroidal construction which allows us to describe the
fibration $\overline{D}_2(n)\to\calA_2^*(n)$ explicitly. Here we shall
concentrate on a description of this map in the most difficult situation,
namely in the neighbourhood of a cusp of maximal corank.

The toroidal compactification $\calA_g^*(n)$ is given by
the second Voronoi decomposition $\Sigma_g$. This is a rational polyhedral
decomposition of the convex hull in $\on{Sym}_g^{\geq0}(\RR)$ of the set
$\on{Sym}_g^{\geq0}(\ZZ)$ of integer semi-positive $(g\times g)$-matrices. For
$g=2$ and $3$ it can be described as follows. First note that $\on{Gl}(g,\ZZ)$
acts on $\on{Sym}_g^{\geq0}(\RR)$ by $\gamma\mapsto
\sideset{^t}{}{\on{\mathit{M}}}\gamma M$. For $g=2$ we define the standard cone
$$
\sigma_2=\RR_{\geq0}\gamma_1+\RR_{\geq0}\gamma_2+\RR_{\geq0}\gamma_3
$$
with
$$
\gamma_1=\begin{pmatrix} 1&0 \\ 0&0 \end{pmatrix},\quad
\gamma_2=\begin{pmatrix} 0&0 \\ 0&1 \end{pmatrix},\quad
\gamma_3=\begin{pmatrix} 1&-1 \\ -1&1 \end{pmatrix}.
$$
Then
$$
\Sigma_2=\{M(\sigma_2);\ M\in\on{Gl}(2,\ZZ)\}.
$$
Similarly for $g=3$ we consider the standard cone
$$
\sigma_3=\RR_{\geq0}\alpha_1+\RR_{\geq0}\alpha_2+\RR_{\geq0}\alpha_3+
\RR_{\geq0}\beta_1+\RR_{\geq0}\beta_2+\RR_{\geq0}\beta_3
$$
with
\begin{gather*}
\alpha_1=\begin{pmatrix} 1&0&0 \\ 0&0&0 \\ 0&0&0 \end{pmatrix},\quad
\alpha_2=\begin{pmatrix} 0&0&0 \\ 0&1&0 \\ 0&0&0 \end{pmatrix},\quad
\alpha_3=\begin{pmatrix} 0&0&0 \\ 0&0&0 \\ 0&0&1 \end{pmatrix},\\
\beta_1=\begin{pmatrix} 0&0&0 \\ 0&1&-1 \\ 0&-1&1 \end{pmatrix},\quad
\beta_2=\begin{pmatrix} 1&0&-1 \\ 0&0&0 \\ -1&0&1 \end{pmatrix},\quad
\beta_3=\begin{pmatrix} 1&-1&0 \\ -1&1&0 \\ 0&0&0 \end{pmatrix}.
\end{gather*}
Then
$$
\Sigma_3=\{M(\sigma_3);\ M\in\on{Gl}(3,\ZZ)\}.
$$
We consider the lattices
$$
N_3 = \ZZ\gamma_1+\ZZ\gamma_2+\ZZ\gamma_3\\
$$
$$
N_6 = \ZZ\alpha_1+\ZZ\alpha_2+\ZZ\alpha_3+\ZZ\beta_1+\ZZ\beta_2+\ZZ\beta_3.
$$
The fans $\Sigma_2$ resp.~$\Sigma_3$ define torus embeddings
$T^3\subset X(\Sigma_2)$ and $T^6\subset X(\Sigma_3)$. We denote the
divisors of $X(\Sigma_3)$ which correspond to the $1$-dimensional simplices
of $\Sigma_3$ by
$\calD^i$. Let $\calD=\calD^1$ be the divisor corresponding to
$\RR_{\geq0}\alpha_3$. An open part of $\calD$ (in the $\CC$-topology) is
mapped to the boundary component $\overline{D}_2(n)$. In order to understand
the structure of $\calD$ we also consider the rank $5$ lattice
$$
N_5=\ZZ\alpha_1+\ZZ\alpha_2+\ZZ\beta_1+\ZZ\beta_2+\ZZ\beta_3\cong
N_6/\ZZ\alpha_3.
$$
The natural projection $\rho:N_{6,\RR}\to N_{5,\RR}$ maps the cones of the fan
$\Sigma_3$ to the cones of a fan $\Sigma_3'\subset N_{5,\RR}$. This fan
defines a torus embedding $T^5=(\calD\setminus
\bigcup\limits_{i\neq1}\calD^i)\subset X(\Sigma_3')=\calD$.

The projection
$$
\begin{array}{rccl}
\lambda:&
N_{6,\RR}\cong\on{Sym}_3(\RR) &\longrightarrow&
N_{3,\RR}\cong\on{Sym}_2(\RR)\\
&\begin{pmatrix} a&b&d \\ b&c&e \\ d&e&f \end{pmatrix} &\longmapsto&
\begin{pmatrix} a&b \\ b&c \end{pmatrix}
\end{array}.
$$
maps $\Sigma_3$ to $\Sigma_2$ and factors through $N_{5,\RR}$. In this
way we obtain an induced map
$$
\begin{array}{ccc}
\calD=X(\Sigma_3') &\longrightarrow& X(\Sigma_2)\\
\cup & & \cup\\
T^5 &\longrightarrow& T^3.
\end{array}
$$
In order to describe this map we first consider the standard simplices
$\sigma_3\subset N_{6,\RR}$ and
$\sigma_2\subset N_{3,\RR}$, resp.~$\sigma_3'=\rho(\sigma_3)\subset N_{5,\RR}$.
On the torus $T^6$ (and similarly on $T^5$ and $T^3$) we introduce coordinates
by
$$
t_{ij}=e^{2\pi i\tau_{ij}/n}\qquad (1\leq i,j\leq 3).
$$
These coordinates correspond to the dual basis
of the basis $U_{ij}^*$ of $\on{Sym}(3,\ZZ)$
where the entries of $U_{ij}^*$ are $1$ in positions $(i,j)$ and $(j,i)$ and
$0$ otherwise. One easily checks that $T_{\sigma_3}\cong\CC^6\subset
X(\Sigma_3)$ and as coordinates on $T_{\sigma_3}$ one can take the coordinates
which correspond to the dual basis of the generators
$\alpha_1,\ldots,\beta_3$. Let us denote these coordinates by
$T_1,\ldots,T_6$. A straightforward calculation shows that the inclusion
$T^6\subset T_{\sigma_3}$ is given by
\begin{equation}\label{formula2.1}
\begin {array}{l@{\qquad}l@{\qquad}l}
T_1 = t_{11}t_{13}t_{12}, & T_2 = t_{22}t_{23}t_{12}, & T_3 =
  t_{33}t_{13}t_{23},\\
T_4 = t_{23}^{-1}, & T_5 = t_{13}^{-1}, & T_6 = t_{12}^{-1}.
\end {array}
\end{equation}
Then $\calD\cap T_{\sigma_3}=\{T_3=0\}$. For genus $2$ the corresponding
embedding $T^3\subset T_{\sigma_2}$ is given by
$$
T_1=t_{11}t_{12},\qquad T_2=t_{22}t_{12},\qquad T_3=t_{12}^{-1}.
$$
Finally we consider $T_{\sigma_3'}\cong\CC^5\subset X(\Sigma_3')$. The
projection $\calD=X(\Sigma_3')\to X(\Sigma_2)$ map $T_{\sigma_3'}$ to
$T_{\sigma_2}$. We can use $T_1,T_2,T_4,T_5,T_6$ as coordinates on
$T_{\sigma_3'}$. Since $\lambda(\alpha_1)=\lambda(\beta_2)=\gamma_1$,
$\lambda(\alpha_2)=\lambda(\beta_1)=\gamma_2$ and $\lambda(\alpha_3)=\gamma_3$
we find that
\begin{equation}\label{formula2.2}
\begin{array}{ccl}
T_{\sigma_3'}\cong\CC^5 &\longrightarrow& T_{\sigma_2}\cong\CC^3\\
(T_1,T_2,T_4,T_5,T_6) &\longmapsto& (T_1T_5,T_2T_4,T_6).
\end{array}
\end{equation}
Given any (maximal dimensional) cone $\sigma'=\rho(\sigma)$ in $\Sigma_3'$ we
can describe the map $T_{\sigma'}\to T_{\lambda(\sigma)}$ in terms of
coordinates by the method described above. In this way we obtain a complete
description of the map $\calD\to X(\Sigma_2)$.

Let us now return to the toroidal compactification $\calA_3^*(n)$ of
$\calA_3(n)$. Let $u_0\subset\QQ^6$ be a maximal isotropic subspace. Then we
obtain the compactification of $\calA_3(n)$ in the direction of the cusp
corresponding to $u_0$ as follows: The parabolic subgroup
$P(u_0)\subset\Gamma_3(n)$ is an extension
$$
1\longrightarrow P'(u_0)\longrightarrow P(u_0)\longrightarrow
P''(u_0)\longrightarrow 1
$$
where $P'(u_0)$ is a lattice of rank $6$. We have an inclusion
$\HH_g/P'(u_0)\subset T^6\subset X(\Sigma_3)$ and we denote the interior of
the closure of $\HH_g/P'(u_0)$ in $X(\Sigma_3)$ by $X(u_0)$. Then $P''(u_0)$
acts on $X(u_0)$ and we obtain a neighbourhood of the cusp corresponding to
$u_0$ by $X(u_0)/P''(u_0)$. We have already described the partial
compactification in the direction of a line (in our case $l_0$). Similarly we
can define a partial compactification in the direction of an isotropic plane
$h_0$. The space $\calA_3^*(n)$ is then obtained by glueing all these partial
compactifications.

The result of Nakamura and Tsushima can then be stated as follows: The
restriction of the map $\pi:\calA_3^*(n)\to\overline{\calA}_3(n)$ to the
boundary component $\overline{D}_2(n)$ admits a factorisation
$$
\diagram
\overline{D}_2(n) \rto^{\pi'} \drto_{\pi} & \calA_2^*(n) \dto^{\pi''}\\
                                          & \overline{\calA}_2(n)
\enddiagram
$$
where $\pi'':\calA_2^*(n)\to\overline{\calA}_2(n)$ is the natural map of the
Voronoi compactification $\calA_2^*(n)$ of $\calA_2(n)$ to the Satake
compactification $\overline{\calA}_2(n)$. The map
$\pi':\overline{D}_2(n)\to\calA_2^*(n)$ is a flat family of
surfaces extending the
universal family over $\calA_2(n)$. In order to describe the fibres over the
boundary points of $\calA_2^*(n)$ recall that every boundary component of
$\calA_2^*(n)$ is isomorphic to the Shioda modular surface $S(n)$. We explain
the \emph{type} of a point $P$ in $\calA_2^*(n)$ as follows:
$$
\begin{array}{lcl}
P\text{ has type I}    &\Longleftrightarrow& P\in\calA_2(n)\\
P\text{ has type II}   &\Longleftrightarrow& P\text{ lies on a smooth fibre
                                             of}\\
                       &&                    \text{a boundary component }S(n)\\
P\text{ has type IIIa} &\Longleftrightarrow& P\text{ is a smooth point on a
                                             singular}\\
                       &&                    \text{fibre of }S(n)\\
P\text{ has type IIIb} &\Longleftrightarrow& P\text{ is a singular point of an
                                             $n$-gon}\\
                       &&                    \text{in }S(n).
\end{array}
$$
Points of type IIIb are also often called \emph{deepest points}.

\begin{proposition}[Nakamura, Tsushima]\label{prop2.4}
Assume $n\geq3$. Let $P$ be a point in $\calA_2^*(n)$ and denote the fibre of
the map $\pi':\overline{D}_2(n)\to\calA_2^*(n)$ over $P$ by $A_P$. Then the
following holds:

\noindent $\on{(i)}$ If $P=[\tau]\in\calA_2(n)$ is of type $\on{I}$ then $A_P$
is a smooth abelian surface, more precisely $A_P\cong A_{n,\tau n}$.

\noindent $\on{(ii)}$ if $P$ is of type $\on{II}$, then $A_P$ is a cycle of
$n$ elliptic ruled surfaces.

\noindent $\on{(iii)}$ If $P$ is of type $\on{IIIa}$, then $A_P$ consists on
$n^2$ copies of $\PP^1\times\PP^1$.

\noindent $\on{(iv)}$ If $P$ is of type $\on{IIIb}$, then $A_P$ consists of
$3n^2$ components. These are $2n^2$ copies of the projective plane $\PP^2$ and
$n^2$ copies of ${\tilde{\PP}}^2$, i.e.~$\PP^2$ blown up in $3$ points in
general position.
\end{proposition}

\begin{Proof}
The proof consists of a careful analysis of the map
$\overline{D}_2(n)\to\calA_2^*(n)$ using the description of the map $\calD\to
X(\Sigma_2)$. For details see \cite[section 4]{Ts}.
\end{Proof}

\begin{urems}
(i) The degenerations of type IIIa and IIIb are usually depicted by the
diagrams
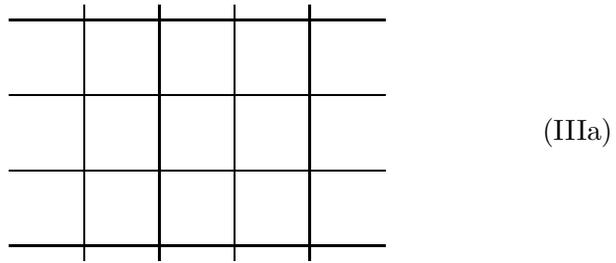
\begin{gather}\tag{IIIa}
\begin{minipage}[c]{3.4cm}
\unitlength1cm
\begin{picture}(5.5,3.4)
\multiput(1,0)(1,0){4}{\line(0,1){3.4}}
\multiput(0,0.2)(0,1){4}{\line(1,0){5}}
\end{picture}
\end{minipage}
\end{gather}
where each square stands for a $\PP^1\times\PP^1$, resp.
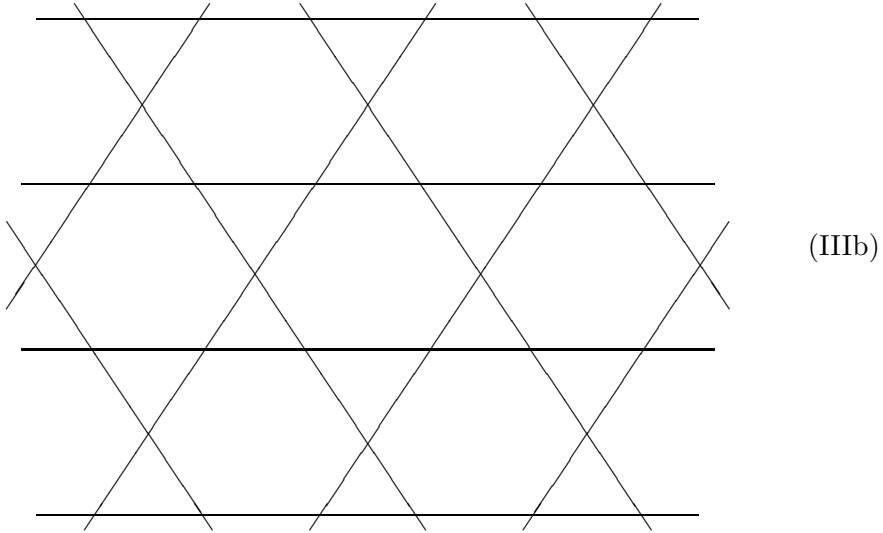
\begin{gather}\tag{IIIb}
\begin{minipage}[c]{9.6cm}
\unitlength1cm
\begin{picture}(9.6,7.5)
\put(0.4,0.2){\line(1,0){8.8}}
\put(0.2,2.4){\line(1,0){9.2}}
\put(0.2,4.6){\line(1,0){9.2}}
\put(0.4,6.8){\line(1,0){8.8}}
\put(0,4.1){\line(2,-3){2.735}}
\put(0.9,7){\line(2,-3){4.67}}
\put(3.9,7){\line(2,-3){4.67}}
\put(6.9,7){\line(2,-3){2.705}}
\put(2.7,7){\line(-2,-3){2.705}}
\put(5.7,7){\line(-2,-3){4.67}}
\put(8.7,7){\line(-2,-3){4.67}}
\put(9.6,4.1){\line(-2,-3){2.735}}
\end{picture}
\end{minipage}
\end{gather}
where the triangles stand for projective planes $\PP^2$ and the hexagons for
blown-up planes ${\tilde{\PP}}^2$.

\noindent (ii) The singular fibres are degenerate abelian surfaces
(cf.~\cite{Na}, \cite{HKW}).

\noindent (iii) This description must be modified for $n=1$ or $2$. Then the
general fibre is a Kummer surface $K_{n,\tau n}$ and the fibres of type (IIIb)
consist of $8$ ($n=2$), resp.~$2$ copies of $\PP^2$.
\end{urems}

The following is a crucial technical step:

\begin{proposition}\label{prop2.5}
Let $n\equiv 0\on{mod} 8p^2$. If $m',m'',\overline{m}',\overline{m}''\in
\frac{1}{2p}\ZZ^2$ then the sections
$\Theta_{m'm''}(\tau,z)\Theta_{\overline{m}'\overline{m}''}(\tau,z)$ of
the line
bundle $M(n)$ on $D_2(n)$ extend to sections of the line bundle
$\overline{M}(n)$ on $\overline{D}_2(n)$.
\end{proposition}

\begin{Proof}
We have to prove that the sections in question extend to the part of
$D_2(n)$ which lies over the boundary of $\calA_2^*(n)$.
This is a local statement. Moreover it is enough to prove extension in
codimension $1$. Due to symmetry considerations we can restrict ourselves to
one boundary component in $\calA_2^*(n)$. We shall use the above description
of the toroidal compactifications $\calA_2^*(n)$ and $\calA_3^*(n)$ and of the
map $\overline{D}_2(n)\to\calA_2^*(n)$. We consider the boundary component of
$\calA_2^*(n)$ given by $\{T_2=0\}\subset T_{\sigma_2}\subset X(\Sigma_2)$.
Recall the theta functions
$$
\Theta_{m'm''}(\tau,z)=\sum_{q\in\ZZ^2} e^{2\pi i[\frac12(q+m'){\tau}
{^t(q+m')}+{(q+m')}{^t(z+m'')}]}
$$
In our situation
$$
\tau=\begin{pmatrix} \tau_{11}&\tau_{12} \\ \tau_{12} & \tau_{22}
\end{pmatrix},\quad z=(z_1,z_2)=(\tau_{13},\tau_{23}).
$$
In level $n$ we have the coordinates
$$
t_{ij}=e^{2\pi i\tau_{ij}/n}
$$
and $\Theta_{m'm''}(\tau,z)$ becomes
\begin{multline*}
\Theta_{m'm''}(\tau,z)=\sum_{q=(q_1,q_2)\in\ZZ^2} t_{11}^{\frac12(q_1+m_1')^2n}
t_{12}^{(q_1+m_1')(q_2+m_2')n} t_{22}^{\frac12(q_2+m_2')^2n}\\
t_{13}^{(q_1+m_1')n} t_{23}^{(q_2+m_2')n} e^{2\pi {i (q+m')} {^tm}''}.
\end{multline*}
We use the coordinates $T_1,T_2,T_4,T_5,T_6$ on $T_{\sigma_3'}$. It follows
from (\ref{formula2.1}) that
\begin{alignat}{2}
t_{11} &= T_1T_5T_6, & \quad t_{22} &= T_2T_4T_6, \notag \\
t_{23} &= T_4^{-1},  & \quad t_{13} &= T_5^{-1},\label{formula2.3}\\
t_{12} &= T_6^{-1}.\notag
\end{alignat}
This leads to the following expression for the theta-functions
\begin{multline*}
\Theta_{m'm''}(\tau,z)=\sum_{q\in\ZZ^2} T_1^{\frac12(q_1+m_1')^2n}
T_2^{\frac12(q_2+m_2')^2n} T_4^{\frac12(q_2+m_2')(q_2+m_2'-2)n}\\
T_5^{\frac12(q_1+m_1')(q_1+m_1'-2)n} T_6^{\frac12((q_1+m_1')-(q_2+m_2'))^2n}
e^{2\pi {i(q+m')} {^tm}''}.
\end{multline*}
By (\ref{formula2.2}) the locus over $T_2=0\subset T_{\sigma_2}$ in
$T_{\sigma_3'}$ is given by $T_2T_4=0$. The equation for the boundary
component $\overline{D}_2(n)$ is given by $t_{33}=0$. Since by
(\ref{formula2.1}) we have $t_{33}=T_3T_4T_5$ we can assume that the normal
bundle and hence $\overline{M}(n)$ (more precisely its pullback to
$X(\Sigma_3')$) is trivial outside $T_4T_5=0$. Since the exponent of $T_2$ is a
non-negative integer (here we use $n\equiv 0\on{mod} 8p^2$)
this shows that the sections extend over $T_2=0$,
$T_4\neq0$. To deal with the other components of $T_{\Sigma_3'}$ which lie
over $\{T_2=0\}$ in $T_{\sigma_2}$ we use the matrices
$$
\nu_{nm}=\begin{pmatrix} 1&0&m \\ 0&1&n \\ 0&0&1 \end{pmatrix}
\qquad (n,m\in\ZZ)
$$
(cf.~\cite{Ts}) which act on $\on{Sym}_3^{\geq0}(\ZZ)$ by
$$
\gamma\longmapsto\sideset{^t}{_{nm}}{\on{\nu}}\gamma\nu_{nm}.
$$
Via $\lambda$ this action lies over the trivial action on
$\on{Sym}_2^{\geq0}(\ZZ)$. This action also factors through $\rho$. Let
$(\sigma_3')_{nm}=\rho(\sideset{^t}{_{nm}}{\on{\nu}}\sigma_3\nu_{nm})$. We
can then either argue with the symmetries induced by this operation or repeat
directly the above calculation for $T_{(\sigma_3')_{nm}}$. Acting with
$\nu_{0m}$, $m\in\ZZ$, we can thus treat all components in $X(\Sigma_3')$
lying over $\{T_2=0\}$ in $X(\Sigma_2)$.
\end{Proof}

\section{Curves in the boundary }

We can now treat curves contained in a boundary component. The
following technical lemma will be crucial.
Its proof uses the
ideas of \cite[Abschnitt~4]{We} in an essential way
and it can be generalized in a suitable form to arbitrary $g$.
We consider the boundary
component $\overline{D}_2(n)$ which belongs to the line
$l_0=(0,\ldots,0,1)\subset\QQ^6$. Recall that the open part $D_2(n)$ of
$\overline{D}_2(n)$ is of the form $D_2(n)=\CC^2\times\HH_2/(P''(l_0)\cap
\Gamma(n))$ and that the group $P''(l_0)/(P''(l_0)\cap\Gamma(n)$) acts on
$\overline{D}_2(n)$. Recall also the fibration $\pi':\overline{D}_2(n)\to
\calA_2^*(n)$. We shall denote the boundary of $\calA_2^*(n)$ by $B$.

\begin{proposition}\label{prop2.6}
Let $(z,\tau)\in\CC^2\times\HH_2$. For every $\eps>0$ there exist integers
$n,k$ and a section $s\in H^0(\overline{M}(n)^k)$ such that

\noindent $\on{(i)}$ $s([z,\tau])\neq0$ where $[z,\tau]\in D_2(n)=\CC^2\times
\HH_2/(P''(l_0)\cap\Gamma(n))$,

\noindent $\on{(ii)}$ $s$ vanishes on $\pi^*B$ of order $\lambda$ with
$\frac{\lambda}{k}\geq\frac{n}{12+\eps}$.
\end{proposition}

\begin{Proof}
Let $p \ge 3$ be a prime number (which will be chosen later).
For $l=2p$
we consider the set of characteristics $\calM$ in $(\frac1l\ZZ/\ZZ)^6$ of the
form $m=(m_p,m_2)$ in $(\frac1p\ZZ/\ZZ)^6\oplus (\frac12\ZZ/\ZZ)^6$ with
$m_p\not\in{\ZZ}^6$. The group $\Gamma_3(1)$ acts on $\calM$ with $2$
orbits. Assume $\eps>0$ is given and that $\widetilde\calM$ is
a subset of $\calM$ with
$$
\#\widetilde\calM<\eps\#\calM.
$$
Then set
$$
\Theta_{\calM,\widetilde\calM}(\tau,z)=
\prod_{m\in\calM\setminus\widetilde\calM}
\Theta_m^l(\tau,z).
$$
Let $n=8p^2$. By Proposition~(\ref{prop2.5})
the functions $\Theta_m^l(\tau,z)$ define sections in
$\overline{M}(n)^p$.
Let $M_1,\ldots,M_N\in
\Gamma_2(1)$ be a set of generators of $\Gamma_2(1)/\Gamma_2(n)\cong
\on{Sp}(4,\ZZ/n\ZZ)$. Then $M_1,\ldots,M_N$, considered as elements in
$P(l_0)$, act on
the line bundle
$\overline{M}(n)$. We set
$$
F_r(\tau,z)=\sum_{i=1}^NM_i^*\Theta_{\calM,\widetilde\calM}^r.
$$
This is a $\Gamma_2/\Gamma_2(n)$-invariant section of $\overline{M}(n)^{pr}$.

Now consider the abelian surface $A=A_{\tau,1}=\CC^2/(\ZZ^2\tau+\ZZ^2)$. Then
$A_{n\tau,n}=\CC^2/((n\ZZ)^2\tau+(n\ZZ)^2)$ is the fibre of $\pi$ over the
point $[\tau]\in\calA_2(n)$. Let
$$
\widetilde\calM=\{m\in\calM;\ \Theta_m(\tau,z)=0\}.
$$
The argument of Weissauer shows that
$$
\#\widetilde\calM<\eps\#\calM
$$
for $p$ sufficiently large. For some $r$ the section $F_r(\tau,z)$
does not vanish at $[z,\tau]\in D_2(n)$.
Let $B'$ be a boundary boundary component of $\calA_2^*(n)$.
The inverse image $D'$ of $B'$ under $\pi'$
consists of several components.
Using the matrices $\nu_{nm}$ which were introduced in the proof of
Proposition~(\ref{prop2.5}) one can, however,
show that the vanishing order of the
sections $\Theta_m^l(\tau,z)$ on the components of $D'$
only depends on $B'$.
Hence one can argue as in \cite{We}
and finds that the vanishing order along $\pi^*B$ goes to
$\frac{prn}{12}$ as $p$ goes
to infinity. Setting $k=pr$ this gives (ii).
\end{Proof}

We can now start giving the proof of Theorem~(\ref{theo0.1}). Let
$$
H=aL-bD\qquad  b>0,\ 12a-\frac bn>0
$$
be a divisor on $\calA_g^*(n)$. In view of Proposition~(\ref{prop1.4}) it
remains to consider curves $C$ which are contained in the boundary.
To simplify
notation we write the decomposition of the boundary $D$ as
$$
D=\sum_{i=1}^N\overline{D}_{g-1}^i(n)
$$
where $N=N(n,g)$ can be computed explicitly. Then
\begin{equation}\label{formula2.4}
H|_{\overline{D}_{g-1}^1(n)}=\left.
\left(aL-b\sum_{i\neq1}\overline{D}_{g-1}^i(n)\right)
\right|_{\overline{D}_{g-1}^1(n)}-
b\overline{D}_{g-1}^{1}(n)|_{\overline{D}_{g-1}^{1}(n)}.
\end{equation}
Now let $g=2$ or $3$ where we have the fibration
$$
\pi':\overline{D}_{g-1}^1(n)\longrightarrow\calA_{g-1}^*(n).
$$
We shall denote the boundary of $\calA_{g-1}^*(n)$ by $B$. Also note that
the restriction of $L$ to the boundary equals ${\pi'}^* L_{\calA_{g-1}^*(n)}$
where we use the notation $L$ for both the line bundle on $\calA_g^*(n)$ and
$\calA_{g-1}^*(n)$. Thus we find that
\begin{equation}\label{formula2.4a}
H|_{\overline{D}_{g-1}^1(n)}=
{\pi'}^*(aL-bB)-b\overline{D}_{g-1}^{1}(n)|_{\overline{D}_{g-1}^{1}(n)}.
\end{equation}
In view of the definition of the line bundle $\overline{M}(n)$ this gives
\begin{equation}\label{formula2.5}
H|_{\overline{D}_{g-1}^1(n)}={\pi'}^*\left(\left(a-\frac
bn\right)L-bB\right)+\frac bn\overline{M}(n).
\end{equation}

\begin{Proof*}{Proof of Theorem~(\ref{theo0.1}) for $g=2$}
In this case the boundary components $\overline{D}_{1}^i(n)$ are
isomorphic to
Shioda's modular surface $S(n)$ and the projection $\pi'$ is just projection
to the modular curve $X(n)$. The degree of $L$ on $X(1)$ is $\frac{1}{12}$ and
we have one cusp. Hence
$$
\deg_{X(n)}(aL-bB)=\mu(n)\left(\frac{a}{12}-\frac bn\right)
$$
where $\mu(n)$ is the degree of the Galois covering $X(n)\to X(1)$,
i.e.~$\mu(n)=|\on{PSL}(2,\ZZ/n\ZZ)|$. This is non-negative if and only if
$a-12\frac bn\geq0$. The normal bundle of
$\overline{D}_{1}^i(n)$ can also be
computed explicitly. This can be done as follows:
Using the degree $10$ cusp form which vanishes on the reducible locus
one finds the equality
$10L=2H_1+D$ on $\calA_2^*$
where $H_1$ is the Humbert surface parametrizing
polarized abelian surfaces which are products.
Hence we conclude for the canonical bundle on $\calA_2^*(n)$
that $K=(3-\frac{10}{n})L+\frac{2}{n}H_1$. The restriction of the divisor
$H_1$ to a boundary component $\overline{D}_{1}^i(n) \cong S(n)$ is the sum of
the
$n^2$  sections $L_{ij}$ of $S(n)$. The canonical bundle of the surfaces $S(n)$
is  equal to the pull-back via ${\pi'}$
of $3L$ minus the divisor of the cusps on the modular
curve $X(n)$ (see also \cite{BH}). Hence adjunction together with an easy
calculation gives
$$
-n\overline{D}_{1}^i(n)|_{\overline{D}_{1}^i(n)}
=2{\pi'}^* L_{X(n)}+2\sum {L_{ij}}
$$
Since
$L_{ij}|_{L_{ij}}=-L_{X(n)}$ one sees immediately that this line bundle is nef
and positive on the fibres of $\pi':S(n)\to X(n)$. The result now follows
directly from (\ref{formula2.4a}).
\end{Proof*}

We shall now turn to the case $g=3$. As we have remarked before
it remains to consider curves which are contained in the
boundary of $\calA_3^*(n)$. Among those curves we shall first deal with curves
whose image under the map $\pi'$ meets the interior of $\calA_2(n)$.

\begin{proposition}\label{prop2.8}
Let $H=aL-bD$ be a divisor on $\calA_3^*(n)$ with $a-12\frac bn>0$, $b>0$. For
every curve $C$ in a boundary component $\overline{D}_2(n)$ with
$\pi'(C)\cap\calA_2(n)\neq\emptyset$ the intersection number $H.C>0$.
\end{proposition}

\begin{Proof}
We shall use (\ref{formula2.5}) and Proposition~(\ref{prop2.6}). If we replace
$n$ by some multiple and consider the pull-back of $H$ the coefficient $b/n$
is not changed. The inverse image of $C$ may have several components.
All of these are, however, equivalent under some finite sympectic group and it
is sufficient to prove that the degree of $H$ is positive on one (and hence
on every) component lying over $C$. After this reduction we can again
assume that $C$ is irreducible and by Proposition~(\ref{prop2.6})
we can find for every $\eps>0$ a divisor ${\cal C}$ not containing $C$ with
$$
\overline{M}(n)={\cal C} +\frac{\lambda}{k} \pi^*B,\qquad \frac{\lambda}{k}\geq
\frac{n}{12+\eps}.
$$
By (\ref{formula2.5})
$$
H|_{\overline{D}_2(n)}=\pi^*\left(\left(a-\frac bn\right)L-b\left(1-
\frac{\lambda}{nk}\right)B\right)+\frac bn {\cal C}.
$$
The assertion follows from the corresponding result for $g=2$ provided
$$
\left(a-\frac bn\right)-12\frac bn\left(1-\frac{\lambda}{nk}\right)\ge
\left(a-12\frac bn\right)-\frac bn \left(1-\frac{12}{12+\eps}\right)>0.
$$
Since $a-12b/n>0$ this is certainly the case for $\eps$ sufficiently small.
\end{Proof}

We are now left with curves in the boundary of $\calA_3^*(n)$ whose
image under $\pi'$ is contained in the boundary of $\calA_2^*(n)$.
These are exactly the curves which are contained in more than $1$
boundary component of $\calA_3^*(n)$.
Before we conclude the proof, we have to analyze the situation once more.
First of all we can assume by symmetry arguments that $C$ is contained in
$\overline{D}_2(n)$=$\overline{D}_{2}^1(n)$.
Let $B'$ be a component
of the boundary $B$ of $\calA_2^*(n)$
which contains $\pi'(C)$. Let
$D'=({\pi'})^{-1}(B')$. Then $D'$ consists of $n$ irreducible components and
we have the following commutative diagram $(n\ge
3)$:

$$
\diagram
  \overline{D}_2^1(n) \rto^{\pi'}
    \morphism{\dottedwith{}}\notip\notip[d]|{\displaystyle\cup}
  &\calA_2^*(n)
    \morphism{\dottedwith{}}\notip\notip[d]|{\displaystyle\cup} \\
  D' \rto^{\pi'} \drto_{\pi}     & B'\cong S(n) \dto^{\pi''} \\
                                 & X(n).
\enddiagram
$$

\noindent Altogether there are three possibilities:\\
(1)
$\pi'(C)=pt$, i.e. $C'$ is contained in a fibre of $\pi'$.\\
(2)
$\pi(C)=pt, \pi'(C) \neq pt$. Then $\pi'(C)$ is either a smooth fibre of $S(n)$
or a component of a singular $n$--gon.\\
(3)
$\pi(C)=X(n)$.\\
The final step in the proof of Theorem (\ref{theo0.1}) is the following:

\begin{proposition}\label{prop2.10}
Let $C\subset \overline{D}_2(n)$ be a curve whose image $\pi'(C)$ is
contained in the
boundary of ${\cal A}_2^{*}(n)$. If $H=aL-bD$ is a divisor with $b>0,
a-12\frac bn>0$ then $H.C>0$.
\end{proposition}

\begin{Proof}
By induction on $g$ and formula (\ref{formula2.4a}) it is enough to prove that
there is some
$\overline{D}_{2}^{j}(n)$ with $C.\overline{D}_{2}^{j}(n)\le 0$. Consider the
inverse image $D'$ of $B'$ under $\pi'$.
Then $D'$ consists of $n$ irreducible components each of
which is of the form
$\overline{D}_2^i(n)\cap\overline{D}_2^1(n)$ for some $i\neq 1$. We already
know
that $-B'|_{B'}$ is nef. Hence
$$
\left(\sum\limits_{i\in I}
\overline{D}_2^i(n)\cap\overline{D}_2^1(n)\right).C\le 0
$$
where $I$ is a suitable set of indices consisting of
$n$ elements. In particular
$\overline{D}_{2}^{j}(n).C\le 0$ for some index $j$.
\end{Proof}
\begin{urems}
(i)
If $\pi'(C)=pt$, then one can
give an alternative proof of $\overline{D}_2(n).C>0$ by computing the
normal bundle of $\overline{D}_2(n)$
restricted to the singular fibres of $\pi'$.
The conormal bundle is ample as in the
smooth case (cf. Lemma (\ref{lemma2.2})). \\
(ii)
If $\pi'(C)\neq pt$ one can also
use the theta functions $\Theta_{m'm''}$ with
$m', m''\in\frac 12 \ZZ^2$ to construct sections of $\overline{M}(n)$ which,
after subtracting suitable components
of the form $\overline{D}_{2}^{i}(n)\cap\overline{D}_2^1(n)$,
do not vanish identically on $C$. In this way one can compute similarly to
the proof of Proposition (\ref{prop2.8}) that
$H.C> 0$.
\end{urems}

\noindent{\em Proof of Theorem } (0.1)(g=3). This follows now immediately from
Proposition (\ref{prop1.4}), Proposition (\ref{prop2.8}) and
Proposition (\ref{prop2.10}). \hfill $\Box$\\

\noindent {\em Proof of the corollaries}. These follow immediately from
Theorem (\ref{theo0.1}) since the moduli spaces are smooth and since
$$
K\equiv(g+1)L-D.
$$
Obviously
$$
(g+1)-\frac{12}n \ge 0\Leftrightarrow\left\{
\begin{array}{ccl}
n \ge 4 &\mbox{ if }& g=2\\
n \ge 3 &\mbox{ if }& g=3.
\end{array}
\right.
$$

\noindent Hence $K$ is nef if $g=2, n\ge 4$ and $g=3, n\ge 3$, resp.
numerically
positive if $g=2, n\ge 5$ and $g=3, n\ge 4$. It follows from general results of
classification theory that $K$ is ample in the latter case.\hfill$\Box$

%
%
\bibliographystyle{alpha}

\vspace{1cm}

\noindent

Author's address:

\bigskip

\parbox{9cm}
{Klaus Hulek\\
Institut f\"ur Mathematik\\
Universit\"at Hannover\\
D 30060 Hannover\\
Germany\\
E-Mail : hulek{\symbol{64}}math.uni-hannover.de}

\end{document}